\newcommand{\al}{\alpha}               
\newcommand{\be}{\beta}
\newcommand{\ga}{\gamma}

\newcommand{\lb}{\lambda}              
\newcommand{\Lb}{\Lambda}

\newcommand{\veps}{\varepsilon}        
\newcommand{\vphi}{\varphi}

\newcommand{\cal}{\mathcal}

\newcommand{\calf}{{\cal F}}

\newcommand{\diam}{{\rm diam}}      

\newcommand{\Fix}{{\rm Fix}}

\newcommand{\incl}{\subseteq}

\newcommand{\es}{\emptyset}          
\newcommand{\sm}{\setminus}
\newcommand{\impl}{\Rightarrow}      
\newcommand{\limpl}{\Longrightarrow}
  
\newcommand{\lequi}{\Longleftrightarrow}
\newcommand{\oo}{\infty}

\newcommand{\sk}{\smallskip}
       
\newcommand{\n}{\noindent}

\def\R+oo{R_+\cup\{\oo\}}

\def\btends   {\stackrel {\it b}{\longrightarrow}}
\def\ctends   {\stackrel {\it c}{\longrightarrow}}
\def\dtends   {\stackrel {\it d}{\longrightarrow}}

\def\gtends   {\stackrel {\it g}{\longrightarrow}}
\def\Dtends   {\stackrel {\it D}{\longrightarrow}}
\def\Gtends   {\stackrel {\it G}{\longrightarrow}}

\newcommand{\barr}{\begin{array}}        
\newcommand{\earr}{\end{array}}
\newcommand{\bcor}{\begin{corollary}}    
\newcommand{\ecor}{\end{corollary}}
\newcommand{\ben}{\begin{enumerate}}     
\newcommand{\een}{\end{enumerate}}
\newcommand{\beq}{\begin{equation}}       
\newcommand{\eeq}{\end{equation}}
\newcommand{\bex}{\begin{example}}        
\newcommand{\eex}{\end{example}}
\newcommand{\bit}{\begin{itemize}}        
\newcommand{\eit}{\end{itemize}}
\newcommand{\blemma}{\begin{lemma}}       
\newcommand{\elemma}{\end{lemma}}
\newcommand{\bproof}{\begin{proof}}       
\newcommand{\eproof}{\end{proof}}
\newcommand{\bprop}{\begin{proposition}}  
\newcommand{\eprop}{\end{proposition}}
\newcommand{\brem}{\begin{remark}}        
\newcommand{\erem}{\end{remark}}
\newcommand{\btab}{\begin{tabular}}       
\newcommand{\etab}{\end{tabular}}
\newcommand{\btheorem}{\begin{theorem}}   
\newcommand{\etheorem}{\end{theorem}}

\documentclass[reqno]{amsart}

\newtheorem{theorem}{\bf Theorem}
\newtheorem{corollary}{\bf Corollary}
\newtheorem{example}{\bf Example}
\newtheorem{lemma}{\bf Lemma}
\newtheorem{proposition}{\bf Proposition}
\newtheorem{remark}{\bf Remark}

\begin{document}

\title
[Contractive Maps in Mustafa-Sims Metric Spaces]
{CONTRACTIVE MAPS IN \\
MUSTAFA-SIMS METRIC SPACES}

\author{Mihai Turinici}
\address{
"A. Myller" Mathematical Seminar;
"A. I. Cuza" University;
700506 Ia\c{s}i, Romania
}
\email{mturi@uaic.ro}




\begin{abstract}
The fixed point result in
Mustafa-Sims metrical structures 
obtained by
Karapinar and Agarwal 
[Fixed Point Th. Appl., 2013, 2013:154]
is deductible from a corresponding one 
stated  in terms of 
anticipative contractions over the associated 
(standard) metric space.
\end{abstract}

\maketitle

\section{Introduction}
\setcounter{equation}{0}

Let $X$ be a nonempty set;
and $d:X\times X\to R_+:=[0,\oo[$
be a {\it metric} over it;
the couple $(X,d)$ is called a {\it metric space}.
Call the subset $Y$ of $X$, 
{\it almost singleton} (in short: {\it asingleton})
provided 
[$y_1,y_2\in Y$ implies $y_1=y_2$];
and {\it singleton},
if, in addition, $Y$ is nonempty;
note that, in this case,
$Y=\{y\}$, for some $y\in X$. 
Further, let $T\in \calf(X)$
be a selfmap of $X$.
[Here, given $A,B\ne \es$, 
$\calf(A,B)$ stands for the class of all functions 
$f:A\to B$; if $A=B$, we write 
$\calf(A)$ in place of  $\calf(A,A)$].
Denote $\Fix(T)=\{x\in X; x=Tx\}$;
each point of this set is referred to as 
{\it fixed} under $T$.
The determination of such points is to be performed
in the context below, comparable with the one in 
Rus \cite[Ch 2, Sect 2.2]{rus-2001}:

{\bf 1a)} 
We say that $T$ is a {\it Picard operator} (modulo $d$) 
if, for each $x\in X$, the iterative sequence
$(T^nx; n\ge 0)$ is $d$-convergent

{\bf 1b)}
We say that $T$ 
is a {\it strong Picard operator} (modulo $d$) 
if, for each $x\in X$,
$(T^nx; n\ge 0)$ is $d$-convergent, 
and $\lim_n(T^nx)$ belongs to $\Fix(T)$

{\bf 1c)} We say that $T$ 
is a {\it globally strong Picard operator} (modulo $d$), 
if it is a strong Picard operator (modulo $d$)
and (in addition), $\Fix(T)$ is an asingleton 
(or, equivalently: singleton).

The sufficient (regularity) conditions 
for such properties are being
founded on {\it orbital} concepts
(in short: o-concepts).
Namely, call the sequence $(z_n; n\ge 0)$ in $X$, 
{\it orbital (modulo $T$)}, when 
it is a subsequence of $(T^n x; n\ge 0)$, 
for some $x\in X$.

{\bf 1d)}
Call $(X,d)$, {\it o-complete}, provided
(for each o-sequence) $d$-Cauchy 
$\limpl$ $d$-convergent 

{\bf 1e)}
We say that $T$ is {\it $(o,d)$-continuous}, if
[$(z_n)$=o-sequence and $z_n\dtends z$] imply  
$Tz_n\dtends Tz$.

When the orbital properties are ignored,
these conventions may be written in the usual way;
we do not give details.

Concerning the existence results for  such points,
a basic one was obtained in  1974 by 
Ciri\'{c} \cite{ciric-1974}.
Call the selfmap $T$, {\it $(d;\al)$-contractive}
(where $\al\ge 0$), when
\ben
\item[] (a01)\ \ 
($\forall x,y\in X$):\
$d(Tx,Ty)\le \\
\al \max\{d(x,Tx),d(x,y),d(x,Ty),d(Tx,y),d(y,Ty)\}$.
\een
Note that, by the very definition of the "max" operator,
it may be also written in the implicit form
\ben
\item[] (a02)\ \ 
($\forall x,y\in X$):\
$d(Tx,Ty)\le \al A(x,y)$;\\
where
$A(x,y)=\diam[T(x;1)\cup T(y;1)]$.
\een
Here, $\diam(U)=\sup\{d(x,y); x,y\in U\}$ 
is the {\it diameter} of the subset $U\incl X$;
and
\ben
\item[]
$T(x;n):=\{T^ix; 0\le i\le n\}$,\ 
$x\in X$, $n\ge 0$;
\een
referred to as: 
the {\it orbital $n$-segment} generated by $x$.

\btheorem \label{t1}
Suppose that $T$ is $(d;\al)$-contractive,
for some $\al\in [0,1[$.
In addition, let $(X,d)$ be o-complete.
Then, $T$ is globally strong Picard (modulo $d$).
\etheorem

This result extends the ones in 
Banach \cite{banach-1922},
Kannan \cite{kannan-1968},
and
Zamfirescu \cite{zamfirescu-1972};
see also 
Hardy and Rogers \cite{hardy-rogers-1973}.
Since all quoted statements have 
a multitude of applications to 
the operator equations theory,
Theorem \ref{t1} was the subject of many extensions.
The most natural one is to pass 
from the "linear" type contraction above
to (implicit) "functional "contractive conditions like
\ben
\item[] (a03)\ \ 
$F(d(Tx,Ty),d(x,y),d(x,Tx),d(y,Ty),d(x,Ty),d(y,Tx))\le 0$,\\ 
for all $x,y\in X$;
\een
where $F:R_+^6\to R$ is a function.
For a basic extension of this type, we refer to 
Dane\v{s} \cite{danes-1976};
further choices of $F$ may be found in 
Rhoades \cite{rhoades-1977}
and the references therein.
Note that, all such conditions
are {\it non-anticipative};
i.e.,
the right member of (a03) does not contain
terms like 
$d(T^iu,T^jv)$, $u,v\in\{x,y\}$, where 
$i+j\ge 3$;
so, the question arises of to what extent 
it is possible to have {\it anticipative} 
contractions (in the above sense).
A positive answer to this 
was recently obtained,
in the "linear" case of Theorem \ref{t1}, 
by 
Dung \cite{dung-2013}.
It is our aim in the present exposition 
to give a further extension of 
this last result, 
in the functional context 
we just quoted.
As an argument for its usefulness,
a fixed point theorem in 
Mustafa-Sims metric spaces
due to 
Karapinar and Agarwal \cite{karapinar-agarwal-2013}
is derived.
This, among others, shows that 
a reduction of their statement to 
standard metrical ones is possible,
along the lines described by 
Jleli and Samet \cite{jleli-samet-2012};
in contradiction with authors' claim.
Further aspects will be delineated elsewhere.

\section{Functional anticipative contractions}
\setcounter{equation}{0}

Let $(X,d)$ be a metric space;
and $T$ be a selfmap of $X$.
In the following, we are interested to 
solve the problem of the introductory part
with the aid of 
(implicit) contractive conditions like
\ben
\item[] (b01)\ \ 
($x,y\in X$):\ $d(Tx,Ty)\le 
\Phi(d(x,Tx),d(x,T^2x),d(x,y),d(x,Ty);\\
d(Tx,T^2x),d(Tx,y),d(Tx,Ty);
d(T^2x,y),d(T^2x,Ty);d(y,Ty))$;
\een
where $\Phi:R_+^{10}\to R_+$ is a certain function.
As precise, these conditions are
anticipative counterparts of 
the (non-anticipative) condition (a03).
To describe them, some conventions are needed.
Given $\vphi\in \calf(R_+)$, we say that $T$ is
{\it anticipative $(d;\vphi)$-contractive}, provided
\ben
\item[] (b02)\ \ 
($\forall x,y\in X$):\
$d(Tx,Ty)\le \vphi(B(x,y))$;\\
where
$B(x,y)=\diam[T(x;2)\cup T(y;1)]$.
\een

The functions $\vphi$ to be considered here 
are to be described as follows.
Call $\vphi\in \calf(R_+)$, {\it increasing}, provided
[$t_1\le t_2$ implies $\vphi(t_1)\le \vphi(t_2)$];
denote the class of all these with $\calf(in)(R_+)$.
The basic properties for such functions to be used in 
the sequel are as follows:

{\bf i)}
Given $\vphi\in \calf(in)(R_+)$, 
we say that it is {\it regressive}, in case
\ben
\item[] 
$\vphi(t)< t$, for all $t> 0$;\ hence, $\vphi(0)=0$.
\een
Note that this property holds in case of
$\vphi$ being {\it super regressive}: 
\ben
\item[]
$\vphi(s+0)< s$, for all $s> 0$;\ hence, $\vphi(0)=0$.
\een
Here, as usually, 
$\vphi(s+0)=\lim_{t\to s+} \vphi(t)$  
is the {\it right limit} of $\vphi$ at $s> 0$.

{\bf ii)}
Call $\vphi\in \calf(in)(R_+)$, 
{\it Matkowski admissible} \cite{matkowski-1975}, provided
\ben
\item[] (b03)\ \ 
$\vphi^n(t)\to 0$ as $n\to \oo$,\ for all  $t> 0$;
\een
here, for each $n\ge 0$, 
$\vphi^n$ stands for the $n$-th {\it iterate} of $\vphi$.
Note that, any such function is regressive;
cf. Matkowski \cite{matkowski-1977}.

{\bf iii)}
For the last one, we need a convention.
Let $\vphi\in \calf(in)(R_+)$ be regressive.
Denote $\psi(t)=t-\vphi(t)$, $t\in R_+$; 
it is an element of $\calf(R_+)$;
referred to as the {\it complement} of $\vphi$.
By definition, the {\it coercive} property for this last function means:
\ben
\item[] (b04)\ \ 
$\lim_{t\to \oo}(\psi(t))=\oo$:\ i.e.:\
$\forall \al> 0$, $\exists \be> \al$:\ 
[$t> \be$ $\limpl$ $\psi(t)> \al$].
\een
By definition, it will be referred to as: 
$\vphi$ is {\it complement coercive}; 
note that, passing to the negation operator, 
this property may be written as:
\ben
\item[] (b05)\ \ 
$\forall \al> 0$, $\exists \be> \al$:\ 
[$t\le \al+\vphi(t)$ $\limpl$ $t\le \be$].
\een

We are now in position to state 
our basic result of this section.

\btheorem \label{t2}
Suppose that $T$ is 
anticipative $(d;\vphi)$-contractive, 
for some regressive, Matkowski admissible, and  
complement-coercive function $\vphi\in \calf(in)(R_+)$.
In addition, let $(X,d)$ be o-complete;
and one of the extra conditions below hold

{\bf 2a)}
$T$ is o-continuous 
[$(x_n)$=o-sequence and $x_n\dtends x$ imply $Tx_n\dtends Tx$]

{\bf 2b)}
$\vphi$ is super regressive.

\n
Then, $T$ is globally strong Picard (modulo $d$); i.e.,

{\bf j)} 
$\Fix(T)=\{z\}$, for some $z\in X$

{\bf jj)}
$T^nx\dtends z$ as $n\to \oo$, for each $x\in X$.
\etheorem

\bproof
We firstly check the asingleton property
of $\Fix(T)$.
Let $z_1,z_1\in \Fix(T)$;
and suppose by contradiction that $z_1\ne z_2$;
hence, $d(z_1,z_2)> 0$.
Clearly,
$$
T(z_1;2)=\{z_1\},\ T(z_2;1)=\{z_2\};\
\mbox{whence},\ B(z_1,z_2)=d(z_1,z_2);
$$
so that, by the contractive condition (and $\vphi$=regressive)
$$
d(z_1,z_2)=d(Tz_1,Tz_2)\le \vphi(d(z_1,z_2))<
d(z_1,z_2);
$$
contradiction. 
Hence, necessarily $z_1=z_2$;
and the asingleton property follows.
It remains now to establish the strong 
Picard property for $T$.
Fix some $x_0\in X$; and put 
$(x_n=T^nx_0$, $n\ge 0)$;
clearly, this is an orbital sequence.
If $x_n=x_{n+1}$ for some $n\ge 0$,
we are done;
so, without loss, one may assume that
\ben
\item[] (b06)\ \ 
$x_n\ne x_{n+1}$ (hence, $\rho_n:=d(x_n,x_{n+1})> 0$), $\forall n$.
\een
Remember that, for each $x\in X$ and each 
$n\ge 0$, 
$T(x;n)=\{T^ix; 0\le i\le n\}$
stands for
the {\it orbital $n$-segment} generated by $x$.  
Put also
\ben
\item[]
$T(x;\oo)=\{T^ix; i\ge 0\}=\cup \{T(x;n); n\ge 0\}$;
\een
and call it: 
the {\it orbital set} generated by $x$.
Note that, by the introduced notations,
we have, for each $k\ge 0$,
\beq \label{201}
T(x_k;n)=\{x_h; 0\le h\le k+n\},\ n\ge 0;\
T(x_k;\oo)=\{x_h; h\ge k\}.
\eeq
Moreover, by the working hypothesis above,
\beq \label{202}
\diam T(x_k; n)\ge \rho_k> 0,\ \mbox{for all}\ k\ge 0,\ n\ge 1.
\eeq
There are several steps to be passed.

{\bf I)}
We start with the following useful evaluation

\blemma \label{le1}
Under the introduced notations,
\beq \label{203}
d(x_i,x_j)\le \vphi(\diam T(x_{i-1};j-i+1)),\ 
\mbox{whenever}\ 1\le i\le j.
\eeq
\elemma

\bproof {\bf (Lemma \ref{le1})}
The case of $i=j$ is clear;
so, without loss, one may assume 
$i< j$; hence, $i+1\le j$.
By definition,
$$  \barr{l}
B(x_{i-1},x_{j-1})=
\diam[T(x_{i-1};2)\cup T(x_{j-1};1)]= 
\diam\{x_{i-1},x_i,x_{i+1},x_{j-1},x_j\} \\
\le \diam\{x_n; i-1\le n\le j\}=\diam T(x_{i-1}; j-i+1);
\earr
$$
wherefrom, combining with the contractive condition,
$$
d(x_i,x_j)\le \vphi(B(x_{i-1},x_{j-1}))\le 
\vphi(\diam T(x_{i-1};j-i+1)).
$$
This ends the argument.
\eproof

{\bf II)}
The following consequence of this fact is to be noted.

\blemma \label{le2}
Put $\al=\rho_0:=d(x_0,x_1) (> 0)$; and let $\be> \al$ be the number 
attached to it, by means of the 
complement-coercive property assumed for $\vphi$. Then,
\beq \label{204}
\diam T(x_0;n)\le \be,\ \ \mbox{for all}\ n\ge 1;
\eeq
hence, necessarily, $\diam T(x_0;\oo)\le \be$.
\elemma

\bproof {\bf (Lemma \ref{le2})}
The case $n=1$ is clear, via $\be> \al=d(x_0,x_1)$;
so, we may assume that $n\ge 2$.
For each $(i,j)$ with $1\le i\le j\le n$, we have, by Lemma \ref{le1},
$$
d(x_i,x_j)\le \vphi(\diam T(x_{i-1};j-i+1))\le
\vphi(\diam T(x_0;n))< \diam T(x_0;n);
$$
so that, necessarily,
$$
\diam T(x_0;n)=d(x_0,x_k),\  \mbox{for some}\  k\in \{1,...,n\}.
$$
On the other hand, the same auxiliary statement gives
$$
d(x_1,x_k)\le \vphi(\diam T(x_0;k))\le \vphi(\diam T(x_0;n)).
$$
Putting these together  yields, by the triangle inequality,
$$  
\barr{l}
\diam T(x_0;n)=d(x_0,x_k)\le d(x_0,x_1)+d(x_1,x_k)\le \\
d(x_0,x_1)+\vphi(\diam T(x_0;n)) =\al+\vphi(\diam T(x_0;n));
\earr
$$
wherefrom, by the complement-coercivity of $\vphi$,
$\diam T(x_0;n)\le \be$; as claimed.
\eproof

{\bf III)}
The following $d$-Cauchy property 
of our iterative sequence is now available.

\blemma \label{le3}
With the same notations as before, one has
\beq \label{205}
\diam T(x_n;\oo)\le \vphi^n(\be),\ \mbox{for all}\ n\ge 0;
\eeq
hence, $(x_n; n\ge 0)$ is a $d$-Cauchy o-sequence.
\elemma

\bproof {\bf (Lemma \ref{le3})}
The case $n=0$ is established in Lemma \ref{le2};
so, we may assume that $n\ge 1$.
By Lemma \ref{le1} one has, for each $(i,j)$ with $n\le i< j$,
$$
d(x_i,x_j)\le \vphi(\diam T(x_{i-1},j-i+1))\le \vphi(\diam T(x_{n-1};\oo)).
$$
Passing to supremum over such $(i,j)$,
yields
$\diam T(x_n;\oo)\le \vphi(\diam T(x_{n-1};\oo))$.
After $n$ steps, one thus gets 
$$
\diam T(x_n;\oo)\le \vphi^n(\diam T(x_0;\oo))\le \vphi^n(\be);
$$
and conclusion follows.
\eproof

{\bf IV)}
As  $(X,d)$ is o-complete,
$x_n\dtends z$, for some 
(uniquely determined) $z\in X$.
There are two alternatives to be discussed.

{\bf Case IV-1.}
Suppose that $T$ is o-continuous.
Then, $(y_n:=Tx_n=x_{n+1}; n\ge 0)$,
$d$-converges to $Tz$.
On the other hand, $(y_n; n\ge 0)$ is a 
subsequence of $(x_n; n\ge 0)$; so that,
$y_n\dtends z$.
As $d$ is sufficient, this yields $z=Tz$.

{\bf Case IV-2.}
Suppose that $\vphi$ is super regressive.
To get the desired fact, we use a 
{\it reductio ad absurdum} argument.
Namely, assume that $z\ne Tz$; i.e.,
$b:=d(z,Tz)> 0$.
From the contractive property, we have
\beq \label{206}
d(x_{n+1},Tz)\le \vphi(B(x_n,z)),
\mbox{for each}\ n\ge 0;
\eeq
where (cf. the previous notations),
$$
B(x_n,z)=\diam[T(x_n,2)\cup T(z;1)]=
\diam\{x_n,x_{n+1},x_{n+2},z,Tz\},\ n\ge 0.
$$
Note that, by the continuity of  the map $(x,y)\mapsto d(x,y)$,
the sequence 
$(\lb_n:=d(x_{n+1},Tz); n\ge 0)$ fulfills $\lb_n\to b$ 
as $n\to \oo$.
On the other hand, by the 
very definition above,
the sequence 
$(\mu_n:=B(x_n,z); n\ge 0)$ fulfills
$$
\mu_n\ge b,\ \forall n;\
\mu_n\to b,\ \mbox{as}\ n\to \oo.
$$
There are two sub-cases to discuss.

{\bf Sub-case IV-2-1.}
Suppose that 
\ben
\item[] (b07)\ \ 
for each  $h\ge 0$, there exists $k> h$, such that $\mu_k=b$.
\een
As a consequence, 
there exists a sequence of ranks $(i(n); n\ge 0)$ 
with $i(n)\to \infty$ as $n\to \infty$, such that  
$\mu_{i(n)}=b$, $\forall n$.
Passing to limit as $n\to \oo$, 
over this subsequence,
in the contractive property (\ref{206}),
yields
$b\le \vphi(b)$; contradiction.

{\bf Sub-case IV-2-2.}
Assume that
the opposite alternative is true:
there exists a certain rank $h\ge 0$, such that
\ben
\item[] (b08)\ \ 
$n\ge h$ $\limpl$ $\mu_n> b$;\ hence
$\mu_n\to b+$ as $n\to \oo$.
\een
Passing to limit in the same
contractive property (\ref{206}), gives 
$b\le \vphi(b+0)< b$; 
again a contradiction.

\n
Summing up, 
the working hypothesis about $z\in X$ 
cannot be accepted;
so, we necessarily have $z=Tz$.
The proof is thereby complete.
\eproof

In particular, assume that the function 
$\vphi$ is linear; i.e.:
$\vphi(t)=\al t$, $t\in R_+$, for some $\al\in [0,1[$.
Then, $\vphi$ is increasing, super regressive, Matkowski 
admissible and complement-coercive.
By Theorem \ref{t2} we get the fixed point statement in 
Dung \cite{dung-2013}.
Given $\al\ge 0$, call $T$, {\it anticipative $(d;\al)$-contractive}, provided
\ben
\item[] (b09)\ \ 
($\forall x,y\in X$):\
$d(Tx,Ty)\le \al B(x,y)$; \\
where (see above)  $B(x,y)=\diam\{x,Tx,T^2x,y,Ty\}$.
\een

\btheorem \label{t3}
Suppose that $T$ is 
anticipative $(d;\al)$-contractive, 
for some $\al\in [0,1[$.
In addition, let $(X,d)$ be o-complete.
Then, $T$ is globally strong Picard (modulo $d$).
\etheorem

{\bf (C)}
For the applications to be considered, 
the following particular case 
of this theorem will be useful.
Denote, for $x,y\in X$,
\ben
\item[] (b10)\ \ 
$P(x,y):=\max\{
d(x,Tx)+d(Tx,y),
d(T^2x,y)+d(T^2x,Ty),\\ 
d(Tx,T^2x)+d(Tx,y), 
d(Tx,y)+d(Tx,Ty),d(x,y),d(x,Ty),d(y,Ty)\}$,
\item[] (b11)\ \ 
$Q(x,y):=\max\{
d(x,Tx)+d(Tx,T^2x), 
d(x,Tx)+d(Tx,y), \\
d(T^2x,Ty)+d(y,Ty), 
d(Tx,T^2x)+d(T^2x,Ty),
d(x,y),d(x,Ty)\}$.
\een
Further, given some $\ga\ge 0$, 
we say that $T$ is 
{\it $(d,P,Q;\ga)$-contractive}, provided
\ben
\item[] (b12)\ \ 
($\forall x,y\in X$):\
$d(Tx,Ty)\le \ga \max\{P(x,y),Q(x,y)\}$.
\een
Note that, by the convention above,
this contractive condition is anticipative. 

The following  fixed  point result is available.

\btheorem \label{t4}
Suppose that $T$ is $(d,P,Q;\ga)$-contractive, for 
some $\ga\in [0,1/2[$.
In addition, let $(X,d)$ be complete.
Then, $T$ is globally strong Picard (modulo $d$).
\etheorem

\bproof
By the very conventions above, one has
$$
P(x,y), Q(x,y)\le 2B(x,y),\ \forall x,y\in X.
$$
So, by the accepted contractive conditions,
it follows that 
$$
(\forall x,y\in X):\ 
d(Tx,Ty)\le 2\ga B(x,y).
$$
Hence, the preceding result applies, with
$\al=2\ga$.
This ends the argument.
\eproof

As a consequence, Theorem \ref{t4} 
is indeed reducible to the developments above.
However, for simplicity reasons,
it would be useful having a 
separate proof of it.

\bproof {\bf (Theorem \ref{t4}) [Alternate]}
First, we establish the 
asingleton property of $\Fix(T)$.
Let $r,s$ be two points in 
$\Fix(T)$. By definition,
$P(r,s)=2d(r,s)$, $Q(r,s)=d(r,s)$;
so that, from the contractive condition,
$$
d(r,s)=d(Tr,Ts)\le 2\ga d(r,s).
$$
This, along with $0\le 2\ga< 1$, 
yields $d(r,s)=0$; whence, $r=s$.
It remains now to establish the 
strong Picard (modulo $d$) property of $T$. 
To this end, we start from 
\beq \label{207}
P(x,Tx),Q(x,Tx)\le d(x,Tx)+d(Tx,T^2x),\ 
\forall x\in X.
\eeq
By the contractive condition, we therefore get
\beq \label{208}
d(Tx,T^2x)\le \be d(x,Tx),\ \forall x\in X;
\eeq
where $0\le \be:=\ga/(1-\ga)< 1$.
Fix some $x_0\in X$; and put 
$(x_n=T^nx_0; n\ge 0)$.
By the above evaluation,
$$
d(x_n,x_{n+1})\le \be^n d(x_0,x_1),\ \forall n.
$$
This tells us that $(x_n; n\ge 0)$ is 
a $d$-Cauchy sequence.
As $(X,d)$ is complete, there must be some
(uniquely determined) 
$r\in X$ such that $x_n\dtends r$.
We claim that $r=Tr$; and this completes the
argument.
By the contractive condition,
\beq \label{209}
d(x_{n+1},Tr)\le \ga\max\{P(x_n,r),Q(x_n,r)\},\
\forall n.
\eeq
But, from the very definitions above,
one has, for all $n\ge 0$,
$$
\barr{l}
P(x_n,r)=\max\{
d(x_n,x_{n+1})+d(x_{n+1},r),
d(x_{n+2},r)+d(x_{n+2},Tr), \\
d(x_{n+1},x_{n+2})+d(x_{n+1},r),
d(x_{n+1},r)+d(x_{n+1},Tr),\\
d(x_{n},r),d(x_{n},Tr),d(r,Tr)\},
\earr
$$
$$
\barr{l}
Q(x_n,r)=\max\{
d(x_n,x_{n+1})+d(x_{n+1},x_{n+2}),
d(x_n,x_{n+1})+d(x_{n+1},r), \\
d(x_{n+2},Tr)+d(r,Tr),
d(x_{n+1},x_{n+2})+d(x_{n+2},Tr), 
d(x_n,r),d(x_n,Tr)\}.
\earr
$$
This yields
$$
\lim_n P(x_n,r)=d(r,Tr), 
\lim_n Q(x_n,r)=2d(r,Tr);
$$
whence, passing to limit in the 
relation (\ref{209}), one gets
$d(r,Tr)\le 2\ga d(r,Tr)$.
As $0\le 2\ga < 1$, this yields 
$d(r,Tr)=0$; so that, $r=Tr$.
The proof is complete.
\eproof

Note that, further extensions of 
the obtained facts are possible,
in the quasi-ordered setting;
we do not give details.
Further aspects may be found in 
Yeh \cite{yeh-1978};
see also 
Popa \cite{popa-2003}.

\section{Dhage metrics}
\setcounter{equation}{0}

As already precise in the introductory part,
there are many generalizations of 
the Banach's fixed point theorem.
Here, we shall be interested in the 
{\it structural} way of extension,
consisting of the "dimensional"
parameters attached
to the ambient metric being increased.
For example, this is the case when 
the initial metric $d(.,.)$ is to be substituted by a 
{\it generalized metric} 
$\Lb:X\times X\times X\to R_+$ 
which fulfills -- at this level --
the conditions imposed to the standard case.
An early construction of this type was 
proposed in 1963 by
Gaehler \cite{gaehler-1963};
the resulting map
$B:X\times X\times X\to R_+$
was referred to as a {\it 2-metric} on $X$.
Short after, this structure was 
intensively used in 
many fixed point theorems,
under the model in
Namdeo et al \cite{namdeo-dubey-tas-2007},
Negoescu \cite{negoescu-1985}
and others; see also
Ashraf \cite[Ch 3]{ashraf-2005},
for a consistent references list.
However, it must be noted that 
this 2-metric 
is not a true generalization of 
an ordinary metric; for --
as shown in
Ha et al \cite{ha-cho-white-1988} --
the associated real function
$B(.,.,.)$ is not $B$-continuous in its arguments.
This, among others, led 
Dhage \cite{dhage-1992}
to construct --
via different geometric reasons --
a new such object.

{\bf (A)}
Let $X$ be some nonempty set.
By a {\it Dhage metric} (in short: D-metric) 
over $X$ we shall mean 
any map $D:X\times X\times X \to R_+$,
with the properties
\ben
\item[] (c01)\ \ 
$D(x,y,z)=D(x,z,y)=D(y,x,z)=D(y,z,x)=$\\
$D(z,x,y)=D(z,y,x)$, $\forall x,y,z\in X$\ \hfill (symmetric)
\item[] (c02)\ \ 
$(x=y=z)$ $\lequi$ $D(x,y,z)=0$ \hfill (reflexive sufficient)
\item[] (c03)\ \ 
$D(x,y,z)\le D(x,y,u)+D(x,u,z)+D(u,y,z)$,\\
for all $x,y,z\in X$ and $u\in X$\ \hfill (tetrahedral).
\een
In this case, the couple $(X,D)$ will be termed a
{\it D-metric} space.

Define a sequential $D$-convergence $(\Dtends)$ 
on $(X,D)$ according to:
$x_n \Dtends x$ iff $D(x_m,x_n,x)\to 0$ as $m,n\to \oo$;
i.e.,
\ben
\item[] (c04)\ \ 
$\forall \veps>0, \exists i(\veps)$:\
$m,n\ge i(\veps)\impl D(x_m,x_n,x)\le \veps$.
\een
Note that this concept obeys the general rules in
Kasahara \cite{kasahara-1976}.
By definition, $x_n \Dtends x$ will be referred to as:
$x$ is the {\it D-limit} of $(x_n)$.
The set of all these will be denoted D-$\lim_n(x_n)$;
if it is nonempty, then $(x_n)$ is called
{\it $D$-convergent};
the class of all $D$-convergent sequences will 
be denoted Conv$(X,D)$.
Further, let the  $D$-Cauchy structure on $(X,D)$
be introduced as:
call the sequence $(x_n)$ in $X$, {\it $D$-Cauchy}, provided
$D(x_m,x_n,x_p)\to 0$ as $m,n,p\to \oo$; i.e.:
\ben
\item[] (c05)\ \ 
$\forall \veps>0, \exists j(\veps)$:
$m,n,p\ge j(\veps) \impl D(x_m,x_n,x_p)\le \veps$.
\een
The class of all these will be indicated as Cauchy$(X,D)$;
it fulfills the general requirements in 
Turinici \cite{turinici-2007}.

By definition, the pair
(Conv$(X,D)$, Cauchy$(X,D)$) will be
called the 
{\it conv-Cauchy structure} attached to $(X,D)$.
Note that, by the properties of $D$,
each $D$-convergent sequence is $D$-Cauchy too;
referred to as: $(X,D)$ is {\it regular}.
The converse is not in general true;
when it holds, we say that 
$(X,D)$ is {\it complete}.
\sk

{\bf (B)}
According to Dhage's
topological results in the area,
this new metric corrects the "bad" 
properties of a 2-metric.
As a consequence, his construction was 
interesting enough so as to be used 
in the deduction of many fixed point results; 
see, for instance,
Dhage \cite{dhage-2000}
and the references therein.
The setting of all these is to be described as below.
Let $(X,D)$ be a 
$D$-metric space; and
$T\in \calf(X)$ be a selfmap of $X$. 
The determination of the points in $\Fix(T)$ 
is to be performed under the lines of Section 1,
adapted to our context:

{\bf 3a)} 
We say that $T$ is a {\it Picard operator} (modulo $D$) 
if, for each $x\in X$, the iterative sequence
$(T^nx; n\ge 0)$ is $D$-convergent

{\bf 3b)}
We say that $T$ 
is a {\it strong Picard operator} (modulo $D$) 
if, for each $x\in X$,
$(T^nx; n\ge 0)$ is $D$-convergent, 
and $D-\lim_n(T^nx)$ belongs to $\Fix(T)$

{\bf 3c)} We say that $T$ 
is a {\it globally strong Picard operator} (modulo $D$), 
if it is a strong Picard operator (modulo $D$)
and (in addition), $\Fix(T)$ is an asingleton 
(or, equivalently: singleton).

Sufficient conditions guaranteeing these properties are of
$D$-metrical type.
The simplest one is the following.
Call $T$, 
{\it $(D;\al)$-contractive} (for some $\al\ge 0$) if
\ben
\item[] (c06)\ \ 
$D(Tx,Ty,Tz)\le \al D(x,y,z)$,\ $\forall x,y,z\in X$.
\een

The following fixed point statement in 
Dhage \cite{dhage-1992}
is the cornerstone of all further developments in the area.

\btheorem \label{t5}
Let $(X,D)$ be complete bounded; and 
$T:X\to X$ be $(D;\al)$-contractive, for some $\al\in [0,1[$.
Then, $T$ is a globally strong Picard operator (modulo $D$).
\etheorem

In the last part of his reasoning, Dhage tacitly used the
$D$-continuity of the application 
$(x,y,z)\mapsto D(x,y,z)$, expressed as
\ben
\item[]
[$x_n\Dtends x$, $y_n\Dtends y$, $z_n\Dtends z$]\
imply $D(x_n,y_n,z_n)\to D(x,y,z)$.
\een
But, as proved in
Naidu, Rao and Rao \cite{naidu-rao-rao-2004},
the described property is not in general valid.
(This must be related with the developments in 
Mustafa and Sims \cite{mustafa-sims-2004},
according to which 
an appropriate construction
of a topological 
and/or uniform structure over 
$(X,D)$ is not in general possible;
we do not give details). 
A conv-Cauchy motivation 
of this negative conclusion
comes from the fact that
the convergence structure
Conv$(X,D)$ attached to 
our $D$-metric space  
is "too large";  
i.e.: for many sequences $(x_n)$ in $X$,
$D-\lim_n(x_n)$ is the whole of $X$.
Returning to the above discussion,
note that --
technically speaking --
it would be possible that 
the conclusion in 
Dhage's fixed point theorem 
be retainable, with a different proof.
However, as results from an
illuminating example provided by
Naidu, Rao and Rao 
\cite{naidu-rao-rao-2005},
this last hope fails as well;
so that, ultimately,
the above stated fixed point result
is not true.
Hence, summing up,
a fixed point theory in $D$-metric spaces 
is not available,
under the admitted conditions 
upon the underlying structure.

\section{Mustafa-Sims metrics}
\setcounter{equation}{0}

The drawbacks of Dhage metrical structures
we just exposed, determined 
Mustafa and Sims \cite{mustafa-sims-2006}
to look for a different perspective upon this matter.
Some basic aspects of it will be described further.
\sk

{\bf (A)}
Let $X$ be a nonempty set.
By a {\it Mustafa-Sims metric} (in short: MS-metric)
on $X$ we mean any map $G:X\times X\times X\to R_+$
with
\ben
\item[] (d01)\ \ 
$G(.,.,.)$ is symmetric and reflexive  
\hfill (see above)
\item[] (d02)\ \ 
$G(x,x,y)=0$ implies $x=y$ \hfill (plane sufficient)
\item[] (d03)\ \ 
$G(x,x,y)\le G(x,y,z)$,\ $\forall x,y,z\in $X, $y\ne z$
\hfill (MS-property)
\item[] (d04)\ \ 
$G(x,y,z)\le G(x,u,u)+G(u,y,z)$,\ $\forall x,y,z,u\in X$
\hfill (MS-triangular).
\een
In this case, the couple $(X,G)$ will be referred to 
as a {\it MS-metric space}.

The following consequences of these axioms are valid.

\bprop \label{p1}
We have, for each $x,y,z,u\in X$, 
\beq \label{401}
G(x,y,z)\le G(x,x,y)+G(x,x,z),\
\eeq
\beq \label{402}
G(x,y,y)\le 2G(x,x,y),\ G(x,x,y)\le 2G(x,y,y)
\eeq
\beq \label{403}
G(x,y,z)\le G(x,u,z)+G(u,y,z)
\eeq
\beq \label{404}
G(x,y,z)\le (2/3)[G(x,u,y)+G(y,u,z)+G(z,u,x)]
\eeq
\beq \label{405}
G(x,y,z)\le G(x,u,u)+G(y,u,u)+G(z,u,u).
\eeq
\eprop

\bproof
{\bf i)}
From (d04) and (d01) we have 
(taking $u=y$),
$G(x,y,z)\le G(x,y,y)+G(z,y,y)$;
this, again via (d01), gives (\ref{401}),
by replacing $(x,y)$ with $(y,x)$.

{\bf ii)}
The first half of (\ref{402}) follows at once 
from (\ref{401}) by taking $z=y$;
and the second part is obtainable 
by replacing $(x,y)$ with $(y,x)$.

{\bf iii)}
By (d01), it results that (d03)
may be written as
\ben
\item[] (d05)\ \ 
$G(x,y,y)\le G(x,y,z)$,\ $\forall x,y,z\in $X, $x\ne z$.
\een
Combining with (d04), we get (for $x\ne z$)
$$ 
G(x,y,z)\le G(x,u,u)+G(u,y,z)\le G(x,u,z)+G(u,y,z);
$$
i.e.: (\ref{403}) holds, when $x\ne z$.
It remains to establish the case $x=z$
of this relation:
\beq \label{406}
G(y,x,x)\le G(x,u,x)+G(y,x,u).
\eeq
Clearly, the alternative $u\ne y$ is 
obtainable from (d05).
On the other hand, the alternative $u=y$ means
$$
G(y,x,x)\le G(x,y,x)+G(y,x,y);
$$
evident. Hence (\ref{403}) is true. 

{\bf iv)}
By a repeated application of (\ref{403}),
$$ \barr{l}
G(x,y,z)\le G(x,u,y)+G(x,u,z), \\
G(x,y,z)\le G(y,u,z)+G(y,u,x),\\ 
G(x,y,z)\le G(z,u,x)+G(z,u,y).
\earr
$$
Adding these, relation (\ref{404}) follows.

{\bf v)}
By (\ref{401}), we have
$$ 
G(y,u,z)\le G(y,u,u)+G(z,u,u).
$$
Replacing in (d04) gives (\ref{405}).
The proof is complete.
\eproof

\brem \label{r1}
\rm

In particular, (\ref{404})
tells us that the MS-metric $G(.,.,.)$ 
is tetrahedral.
Moreover, $G(.,.,.)$ is sufficient.
In fact, assume that 
$G(x,y,z)=0$; but, e.g., $y\ne z$.
From (d03), we then get $G(x,x,y)=0$;
wherefrom, by (d02), $x=y$.
In this case, 
the working hypothesis becomes
$G(y,y,z)=0$; so, again via (d02), $y=z$;
contradiction.
Hence, summing up, 
$G(.,.,.)$ is a D-metric on $X$.
\erem

{\bf (B)}
By an {\it almost metric} on $X$, we mean any map
$g:X\times X\to R_+$ with 
\ben
\item[] (d06)\ \ 
$g(x,y)\le g(x,z)+g(z,y)$, $\forall x,y,z\in X$\
\hfill ({\it triangular})
\item[] (d07)\ \ 
$x=y$ $\lequi$ $g(x,y)=0$ 
\hfill ({\it reflexive sufficient});
\een
see also
Turinici \cite{turinici-2011}.
Some basic examples of such objects 
are to be obtained, 
in the MS-metric space $(X,G)$, as follows.
Define a quadruple of maps $b,c,d,e:X\times X\to R_+$ 
according to:
for each $x,y\in X$,
\ben
\item[] (d08)\ \ 
$b(x,y)=G(x,y,y)$,\ $c(x,y)=G(x,x,y)=b(y,x)$
\item[] (d09)\ \ 
$d(x,y)=\max\{b(x,y),c(x,y)\}$,\ 
$e(x,y)=b(x,y)+c(x,y)$.
\een

\bprop \label{p2}
Under the above conventions,

{\bf j)}
The mappings $b(.,.)$ and $c(.,.)$ are
triangular and reflexive sufficient; 
hence, these are almost metrics on $X$

{\bf jj)}
The mappings $d(.,.)$ and $e(.,.)$ are 
triangular, reflexive sufficient and
symmetric; hence, these are
(standard) metrics on $X$

{\bf jjj)}
In addition, the following relations are valid
\beq \label{407}
b\le 2c\le 2d\le 4b, c\le 2b\le 2d\le 4c,\ 
d\le e\le 2d.
\eeq
\eprop

\bproof
{\bf j)}
It will suffice establishing the 
assertions concerning the map $b(.,.)$.
The reflexive sufficient property is a direct
consequence of  (d01) and (d02).
On the other hand, the triangular property
is a direct consequence of (d04).
In fact, by this condition, we have
(taking $y=z$)
$$
G(x,y,y)\le G(x,u,u)+G(u,y,y);
$$
and, from this we are done.

{\bf jj)}
Evident, by the involved definition.

{\bf jjj)}
The first and second part are immediate, by 
Proposition \ref{p1}.
The third part is evident.
Hence the conclusion.
\eproof

\brem \label{r2}
\rm

A formal verification of 
{\bf j)} 
is to be found in 
Jleli and Samet \cite{jleli-samet-2012}.
On the other hand, {\bf jj)} (modulo $e$)
was explicitly asserted in
Mustafa and Sims \cite{mustafa-sims-2006}.
This determines us to conclude that 
{\bf j)} is also clarified 
by the quoted authors.
\erem

{\bf (C)}
Having these precise, we may now pass to 
the conv-Cauchy structure of 
a MS-metric space $(X,G)$.

Define a sequential $G$-convergence $(\Gtends)$ 
on $(X,M)$ according to:
$x_n \Gtends x$ iff $G(x_m,x_n,x)\to 0$ 
as $m,n\to \oo$; i.e.,
\ben
\item[] (d10)\ \ 
$\forall \veps>0, \exists i(\veps)$:\
$m,n\ge i(\veps) \impl G(x_m,x_n,x)\le \veps$.
\een
As before, this concept obeys the general rules in
Kasahara \cite{kasahara-1976}.
By definition, $x_n \Gtends x$ will be referred to as:
$x$ is the {\it G-limit} of $(x_n)$.
The set of all these will be denoted G-$\lim_n(x_n)$;
if it is nonempty, then $(x_n)$ is called
{\it $G$-convergent};
the class of all $G$-convergent sequences will 
be denoted Conv$(X,G)$.
Call the convergence $(\Gtends)$, {\it separated}
when G-$\lim_n(x_n)$ is an asingleton, for each
sequence $(x_n)$ of $X$.
Further, let the  $G$-Cauchy structure on $(X,G)$
be introduced as:
call $(x_n)$, {\it $G$-Cauchy}, provided
$G(x_m,x_n,x_p)\to 0$ as $m,n,p\to \oo$; i.e.:
\ben
\item[] (d11)\ \ 
$\forall \veps>0, \exists j(\veps)$:
$m,n,p\ge j(\veps) \impl G(x_m,x_n,x_p)\le \veps$.
\een
The class of all these will be indicated as Cauchy$(X,G)$;
it fulfills the general requirements in 
Turinici \cite{turinici-2007}.
By definition, the pair
(Conv$(X,G)$, Cauchy$(X,G)$) will be
called the 
{\it conv-Cauchy structure} attached to $(X,G)$.
Call $(X,G)$, {\it regular} when 
each $G$-convergent sequence is $G$-Cauchy too;
and {\it complete}, if the converse holds: 
each $G$-Cauchy sequence is $G$-convergent.

In parallel to this, we may introduce a 
conv-Cauchy structure attached to any $g\in \{b,c,d,e\}$.
This, essentially, consists in the following.
Define a sequential $g$-convergence $(\gtends)$ 
on $(X,g)$ according to:
$x_n \gtends x$ iff $g(x_n,x)\to 0$. 
This will be referred to as:
$x$ is the {\it g-limit} of $(x_n)$.
The set of all these will be denoted g-$\lim_n(x_n)$;
if it is nonempty, then $(x_n)$ is called
{\it $g$-convergent};
the class of all $g$-convergent sequences will 
be denoted Conv$(X,g)$.
Call the convergence $(\gtends)$, {\it separated}
when g-$\lim_n(x_n)$ is an asingleton, for each
sequence $(x_n)$ of $X$.
Further, let the  $g$-Cauchy structure on $(X,g)$
be introduced as:
call the sequence $(x_n)$ in $X$, {\it $g$-Cauchy}, provided
$g(x_m,x_n)\to 0$ as $m,n\to \oo$;
the class of all these will be indicated as Cauchy$(X,g)$.
By definition, (Conv$(X,g)$, Cauchy$(X,g)$) will be
called the conv-Cauchy structure attached to $(X,g)$.
Call $(X,g)$, {\it regular}, when 
each $g$-convergent sequence is $g$-Cauchy;
and {\it complete}, when the converse holds:
each $g$-Cauchy sequence is $g$-convergent.

\bprop \label{p3}
Under the above conventions,

{\bf i)}
($\forall (x_n)\incl X$, $\forall x\in X$):\
$x_n\Gtends x$ is equivalent with 
\beq \label{408}
\mbox{
$x_n\gtends x$\ 
for some/all $g\in \{b,c,d,e\}$
}
\eeq

{\bf ii)}
the convergence structures $(\Gtends)$ and 
$(\gtends)$ (for $g\in \{b,c,d,e\}$) 
are separated.
\eprop

\bproof
{\bf i)}
Assume that $x_n\Gtends x$;
i.e.: $G(x_m,x_n,x)\to 0$ as $m,n\to \oo$.
This yields
$G(x_n,x_n,x)\to 0$ as $n\to \oo$;
i.e.: $x_n\ctends x$;
wherefrom, combining with Proposition \ref{p2},
(\ref{408}) is clear.
Conversely, assume that (\ref{408}) holds.
Taking  (\ref{407}) into account, gives 
$x_n\btends x$;
wherefrom, by means of (\ref{401}),
$x_n\Gtends x$.

{\bf ii)}
Clearly, $(\gtends)$ is separated, 
for $g\in \{d,e\}$.
This, by the preceding step, gives the 
desired fact.
\eproof

Likewise, the following characterization of the
Cauchy property is available.

\bprop \label{p4}
The following are valid: 

{\bf j)}
($\forall (x_n)\incl X$):\
$(x_n)$ is $G$-Cauchy is equivalent with
\beq \label{409}
\mbox{
$(x_n)$ is $g$-Cauchy, for some/all $g\in \{b,c,d,e\}$
}
\eeq

{\bf jj)}
$(X,G)$ and $(X,g)$ (for $g\in \{b,c,d,e\}$) are 
regular

{\bf jjj)}
$(X,G)$ is complete iff $(X,g)$ is complete, for 
some/all $g\in \{b,c,d,e\}$.
\eprop

\bproof
{\bf j)}
Assume that $(x_n)$ is G-Cauchy;
i.e.: $G(x_m,x_n,x_p)\to 0$ as $m,n,p\to \oo$.
This, in particular, yields
$G(x_m,x_n,x_n)\to 0$ as $m,n\to \oo$;
i.e.: $(x_n)$ is $b$-Cauchy;
so that, combining with Proposition \ref{p2},
(\ref{409}) follows.
Conversely, assume that (\ref{409}) holds.
Taking (\ref{407}) into account, 
one gets  that $(x_n)$ is $c$-Cauchy;
wherefrom, by means of (\ref{401}),
we are done.

{\bf jj)}
The assertion is clear for $(X,G)$, 
by Proposition \ref{p1};
as well as for $(X,g)$ (where $g\in \{d,e\}$),
by its metric properties.
The remaining situations ($g\in \{b,c\}$)
follow from Proposition \ref{p3} and {\bf j)} above.

{\bf jjj)}
Evident, by the previously obtained facts.
\eproof

{\bf (D)}
Let $(X,G)$ be a MS-metric space.
Given the function $\Lb:X\times X\times X\to R$, 
call it {\it sequentially $G$-continuous}, provided
\ben
\item[]
$x_n\Gtends x$, $y_n\Gtends y$, $z_n\Gtends z$
imply $\Lb(x_n,y_n,z_n)\to \Lb(x,y,z)$.
\een
A basic example of this type is 
just the one of $G(.,.,.)$.
To verify this, the following
auxiliary fact is to be used
(cf. Mustafa and Sims \cite{mustafa-sims-2006}):

\bprop \label{p5}
The map $G(.,.,.)$ is $d$-Lipschitz:
\beq \label{410}
|G(x,y,z)-G(u,v,w)|\le d(x,u)+d(y,v)+d(z,w),\
\forall x,y,z,u,v,w\in X;
\eeq
hence, in particular, $G(.,.,.)$ is $d$-continuous.
\eprop

\bproof
By the MS-triangular property of $G$,
$$ \barr{l}
G(u,v,w)\le G(v,y,y)+G(y,u,w),\\
G(u,w,y)\le G(u,x,x)+G(x,y,w),\\
G(w,x,y)\le G(w,z,z)+G(z,x,y);
\earr
$$
so that (by the adopted notations) 
$$
G(u,v,w)-G(x,y,z)\le d(u,x)+d(v,y)+d(w,z).
$$
In a similar way, one gets 
(by replacing $(x,y,z)$ with $(u,v,w)$)
$$
G(x,y,z)-G(u,v,w)\le d(x,u)+d(y,v)+d(z,w).
$$
These, by the symmetry of $d(.,.)$, 
give the written conclusion.
\eproof

As a direct consequence of this, we have
(taking Proposition \ref{p3} into account)

\bprop \label{p6}
The map $G(.,.,.)$ is 
sequentially $G$-continuous in its variables.
\eprop

This property allows us to 
get a partial answer to a useful global question.
Call the MS-metric $G(.,.,.)$, {\it symmetric} if
\ben
\item[]
$G(x,y,y)=G(x,x,y)$,\ $\forall x,y\in X$.
\een
Note that, under the conventions above, this may be
expressed as: $b=c$; wherefrom: $d=b=c$, $e=2b=2c$. 
The class of symmetric MS-metrics is nonempty.
For example, given the metric $g(.,.)$ on $X$, 
its associated MS-metric
\ben
\item[]
$G(x,y,z)=\max\{g(x,y),g(y,z),g(z,x)\}$,\ 
$x,y,z\in X$
\een
is symmetric, as it can be directly seen.
On the other hand, the class of all non-symmetric 
MS-metrics is also nonempty;
see
Mustafa and Sims \cite{mustafa-sims-2006}
for an appropriate example.
Hence, the question of a certain MS-metric on $X$
being or not symmetric is not trivial.
An appropriate answer to this may be given 
as follows.
Call the MS-metric space $(X,G)$, 
{\it perfect} provided 
\ben
\item[]
for each $x\in X$ there exists a sequence
$(x_n)$ in $X\sm \{x\}$ with $x_n\Gtends x$.
\een

\bprop \label{p7}
Suppose that $(X,G)$ is perfect.
Then, $G(.,.,.)$ is symmetric.
\eprop

\bproof
Let $x,y\in X$ be arbitrary fixed.
Further, let $(y_n)$ be a sequence in 
$X\sm \{y\}$ with $y_n\Gtends y$.
From the MS-property of $G(.,.,.)$,
$$
\mbox{
$G(x,x,y)\le G(x,y,y_n)$,\ for all $n$.
}
$$
Passing to limit as $n\to \oo$ yields, via
Proposition \ref{p6}, $G(x,x,y)\le G(x,y,y)$.
As $x,y\in X$ were arbitrary, one gets
(under our notations) 
$c(x,y)\le b(x,y)$, $\forall x,y\in X$.
This gives 
(by a substitution of $(x,y)$ with $(y,x)$),
$b(x,y)\le c(c,y)$; wherefrom $b=c$.
The proof is complete.
\eproof

It follows from this that the class of all
non-symmetric MS-metrics over $X$ is not very large.
Further aspects will be delineated elsewhere.

\section{Contractions over MS-metric spaces}
\setcounter{equation}{0}

From the developments above, it follows that
the metrical construction proposed by 
Mustafa and Sims \cite{mustafa-sims-2006}
corrects certain errors of 
the preceding one, due to
Dhage \cite{dhage-1992}.
As a consequence,
this structure was 
intensively used in 
many fixed point theorems;
see, for instance,
Aage and Salunke \cite{aage-salunke-2012},
Aydi, Shatanawi and Vetro 
\cite{aydi-shatanawi-vetro-2011},
Choudhury and Maity 
\cite{choudhury-maity-2011},
Mustafa, Obiedat and Awawdeh
\cite{mustafa-obiedat-awawdeh-2008},
Saadati et al 
\cite{saadati-vaezpour-vetro-rhoades-2010},
Shatanawi \cite{shatanawi-2010},
to quote only a few.
But, as recently proved by 
Jleli and Samet \cite{jleli-samet-2012},
most fixed point results on MS-metric spaces
are directly reducible to their corresponding 
statements on almost/standard metric spaces.
Under this perspective, the return 
of many writers in the area
to the initial setting of 2-metric spaces 
must be not surprising;
see, for instance,
Dung et al \cite{dung-hieu-ly-thinh-2013}.
Clearly, it is possible that 
{\it not} all fixed point results 
over MS-metric spaces be obtainable in this way;
to verify the claim,
an example was proposed by
Karapinar and Agarwal \cite{karapinar-agarwal-2013}.
Some conventions are in order.
Let $(X,G)$ 
be a MS-metric space; with, in addition,
\ben
\item[] (e01)\ \ 
$(X,G)$ is complete;\ hence, so is $(X,d)$.
\een
Here, $d$ is the associated to $G$
standard metric we just introduced; namely,
\beq \label{501} 
d(x,y)=\max\{G(x,y,y),G(x,x,y)\},\ x,y\in X.
\eeq
Remember that, by the MS-triangular inequality,
$$
G(x,y,z)\le G(x,y,y)+G(y,y,z),\ \forall x,y,z\in X;
$$
and this gives
the so-called {\it strong triangle inequality}:
\beq \label{502}
G(x,y,z)\le d(x,y)+d(y,z),\ \forall x,y,z\in X.
\eeq
Further, let $T$ be a selfmap of $X$.
The question of determining its fixed points
is to be treated with the aid of 
Picard concepts in Section 3 (modulo $G$). 
Sufficient conditions for these properties are 
$G$-counterparts of the ones in Section 1.

{\bf (A)}
Define, for $x,y,z\in X$,
\ben
\item[] (e02)\ \ 
$M(x,y,z)=\max\{
G(x,Tx,y),G(y,T^2x,Ty),G(Tx,T^2x,Ty), \\
G(y,Tx,Ty),G(x,Tx,z),G(z,T^2x,Tz), \\
G(Tx,T^2x,Tz),G(z,Tx,Ty),G(x,y,z), \\
G(x,Tx,Tx),G(y,Ty,Ty),G(z,Tz,Tz), \\
G(z,Tx,Tx),G(x,Ty,Ty),G(y,Tz,Tz)\}$.
\een
Given $\ga\ge 0$, let us say that $T$ is
{\it $(G,M;\ga)$-contractive}, provided
\ben
\item[] (e03)\ \ 
$G(Tx,Ty,Tz)\le \ga M(x,y,z)$,\ 
$\forall x,y,z\in X$.
\een

We are now in position to state the
announced result (in our notations):

\btheorem \label{t6}
Suppose that $T$ is $(G,M;\ga)$-contractive, for some
$\ga\in [0,1/2[$.
Then,

{\bf i)}
$T$ is a globally strong Picard operator (modulo $d$)

{\bf ii)}
$T$ is a globally strong Picard operator (modulo $G$).
\etheorem

According to the authors, 
this fixed point statement is 
an illustration of the following assertion:
there are {\it many} fixed point results over
MS-metric structures,
to which the reduction techniques in 
Jleli and Samet (see above)
are not applicable.
It is our aim in the following to show that, 
actually, the above stated fixed point theorem 
cannot be viewed as such an exception
[i.e.: as an illustration of this 
(hypothetical for the moment) alternative].
Precisely, we shall establish 
that Theorem \ref{t6} is reducible to the
anticipative fixed point result over
standard metric spaces given in a preceding place.
This will follow from the proposed

\bproof {\bf (Theorem \ref{t6})} 
By the accepted condition,
\beq \label{503}
G(Tx,Ty,Ty)\le \ga M(x,y,y),\
G(Tx,Tx,Ty)\le \ga M(x,x,y),\ \forall x,y\in X;
\eeq
and this, by the definition of $d$ (see above), gives
\beq \label{504}
d(Tx,Ty)\le \ga \max\{M(x,y,y),M(x,x,y)\},\
\forall x,y\in X.
\eeq
Now, let us deduce from this 
contractive relation in terms of $G$, a corresponding
contractive relation in terms of $d$.
To this end, we have 
$$
\barr{l}
M(x,y,y)=\max\{
G(x,Tx,y),G(y,T^2x,Ty),G(Tx,T^2x,Ty), \\
G(y,Tx,Ty),G(x,Tx,y),G(y,T^2x,Ty), \\
G(Tx,T^2x,Ty),G(y,Tx,Ty),G(x,y,y), \\
G(x,Tx,Tx),G(y,Ty,Ty),G(y,Ty,Ty), \\
G(y,Tx,Tx),G(x,Ty,Ty),G(y,Ty,Ty)\};
\earr
$$
or equivalently (eliminating the identical terms)
$$
\barr{l}
M(x,y,y)=\max\{
G(x,Tx,y),G(y,T^2x,Ty),G(Tx,T^2x,Ty), \\
G(y,Tx,Ty),G(x,y,y),G(x,Tx,Tx),\\
G(y,Ty,Ty),G(y,Tx,Tx),G(x,Ty,Ty)\}.
\earr
$$
By the strong triangle inequality,
$$
\barr{l}
G(x,Tx,y)\le d(x,Tx)+d(Tx,y), \\
G(y,T^2x,Ty)\le d(T^2x,y)+d(T^2x,Ty), \\
G(Tx,T^2x,Ty)\le d(Tx,T^2x)+d(T^2x,Ty), \\
G(y,Tx,Ty)\le d(Tx,y)+d(Tx,Ty);
\earr
$$
and this yields (by avoiding the smaller terms)
\beq \label{505}
\barr{l}
M(x,y,y)\le P(x,y):=\max\{
d(x,Tx)+d(Tx,y), \\
d(T^2x,y)+d(T^2x,Ty), 
d(Tx,T^2x)+d(Tx,y), \\
d(Tx,y)+d(Tx,Ty),d(x,y),d(x,Ty),d(y,Ty)\}.
\earr
\eeq
Similarly, we have
$$
\barr{l}
M(x,x,y)=\max\{
G(x,Tx,x),G(x,T^2x,Tx),G(Tx,T^2x,Tx), \\
G(x,Tx,Tx),G(x,Tx,y),G(y,T^2x,Ty), \\
G(Tx,T^2x,Ty),G(y,Tx,Tx),G(x,x,y), \\
G(x,Tx,Tx),G(x,Tx,Tx),G(y,Ty,Ty), \\
G(y,Tx,Tx),G(x,Tx,Tx),G(x,Ty,Ty)\};
\earr
$$
or equivalently (eliminating the identical terms)
$$
\barr{l}
M(x,x,y)=\max\{
G(x,x,Tx),G(x,Tx,T^2x),G(Tx,Tx,T^2x), \\
G(x,Tx,Tx),G(x,Tx,y),G(y,T^2x,Ty),\\
G(Tx,T^2x,Ty),G(y,Tx,Tx),G(x,x,y),\\
G(y,Ty,Ty),G(x,Ty,Ty)\}.
\earr
$$
By the strong triangle inequality
$$
\barr{l}
G(x,Tx,T^2x)\le d(x,Tx)+d(Tx,T^2x), \\
G(x,Tx,y)\le d(x,Tx)+d(Tx,y), \\
G(y,T^2x,Ty)\le d(T^2x,Ty)+d(y,Ty), \\
G(Tx,T^2x,Ty)\le d(Tx,T^2x)+d(T^2x,Ty);
\earr
$$
and this yields (by avoiding the smaller terms)
\beq \label{506}
\barr{l}
M(x,x,y)\le Q(x,y):=\max\{
d(x,Tx)+d(Tx,T^2x), \\
d(x,Tx)+d(Tx,y), 
d(T^2x,Ty)+d(y,Ty), \\
d(Tx,T^2x)+d(T^2x,Ty),d(x,y),d(x,Ty)\}.
\earr
\eeq
Summing up, we therefore have
\beq \label{507}
d(Tx,Ty)\le \ga \max\{P(x,y),Q(x,y)\},\
\forall x,y\in X.
\eeq
In other words, 
$T$ is $(d,P,Q;\ga)$-contractive
(according to a preceding convention).
But then,
the metrical fixed point result 
(involving anticipative contractions)
we just evoked
gives us the conclusion
in terms of $d$.
The remaining conclusion 
(in terms of $G$) 
is a direct consequence of it,
by the properties of the Mustafa-Sims convergence we
already sketched.
\eproof

Note, finally, that this reduction process 
comprises as well another fixed point result over
Mustafa-Sims metric spaces
given by 
Karapinar and Agarwal \cite{karapinar-agarwal-2013};
we do not give details.
Further aspects may be found in
Samet et al \cite{samet-vetro-vetro-2012}.


\end{document}